\documentclass{amsart}
\usepackage[latin1]{inputenc}
\newtheorem{thm}{Theorem}
\newtheorem{cor}[thm]{Corollary}
\newtheorem{lem}[thm]{Lemma}
\newtheorem{prop}[thm]{Proposition}
\theoremstyle{definition}
\newtheorem{defn}[thm]{Definition}
\theoremstyle{remark}

\numberwithin{equation}{section}

\newcommand{\To}{\longrightarrow}

\begin{document}

\setcounter{tocdepth}{1}
\thanks{Author supported by FPU grant and grant MTM2005-08379 from
MEC of Spain and grant 00690/PI/2004 of Fundación Séneca of Región
de Murcia} 
\title[]{Linearly ordered Radon-Nikod\'ym compact spaces}%
\author{Antonio Avil\'es}%
\address{Antonio Avilés\\ Departamento de Matemáticas\\ Universidad de Murcia\\ 30100 Espinardo (Murcia)\\ Spain }%
\email{avileslo@um.es}%

\begin{abstract}
We prove that every fragmentable linearly ordered compact space is
almost totally disconnected. This combined with a result of
Arvanitakis yields that every linearly ordered quasi
Radon-Nikod\'ym compact space is Radon-Nikod\'ym, providing a new
partial answer to the problem of continuous images of
Radon-Nikod\'ym
compacta.
\end{abstract}

\maketitle

It is an open problem posed by Namioka~\cite{Namioka} whether the
class of Radon-Nikod\'ym compact spaces is closed under continuous
images. Several authors~\cite{FabHeiMat} \cite{Arvanitakis}
\cite[p. 104]{ArkhangelskiiGT2} who have studied this problem have
introduced some superclasses of the class of Radon-Nikod\'ym
compacta which are closed under continuous images, although all
these classes turned out to be equal to the class of quasi
Radon-Nikod\'ym compacta as shown in~\cite{Namiokanote}
and~\cite{cardinalb}. Let us recall that\\

\begin{enumerate}
\item A compact space $K$ is \emph{Radon-Nikod\'ym compact} if and
only if there exists a lower semicontinuous metric $d:K\times K
\To [0,+\infty)$ which fragments $K$.\\

\item A compact space $K$ is \emph{quasi Radon-Nikod\'ym compact}
if and only if there exists a lower semicontinuous quasi metric
$d:K\times K \To [0,+\infty)$ which fragments $K$.\\

\item A compact space $K$ is a \emph{fragmentable compact} if and
only if there exists a quasi metric $d:K\times K \To [0,+\infty)$
which fragments $K$.\\

\end{enumerate}

Here, a quasi metric is a symmetric map $d:K\times K\To
[0,+\infty)$ such that $d(x,y)=0$ if and only if $x=y$ but which
may fail triangle inequality. Also, a map $d:K\times K \To
[0,+\infty)$ is said to fragment the topological space $K$ if for
every nonempty (closed) subset $L$ of $K$ and every
$\varepsilon>0$ there exists a relative open subset $U$ of $L$ of
diameter less than $\varepsilon$, that is, $\sup\{d(x,y) : x,y\in
U\}<\varepsilon$. Lower semicontinuity means that
the set $\{(x,y) : d(x,y)\leq a\}$ is closed for every $a\geq 0$.\\

The class of fragmentable compacta is larger than the other two,
for instance any Gul'ko non Eberlein compact is an example of
fragmentable and not quasi Radon-Nikod\'ym compact. It is again an
open problem whether every quasi Radon-Nikod\'ym compact is
Radon-Nidkod\'ym compact (as mentioned earlier, the class of quasi
Radon-Nikod\'ym compacta is closed under continuous images, and it
is even unknown whether every quasi Radon-Nikod\'ym compact is the
continuous image of a Radon-Nikod\'ym compact). Mainly two partial
answers to this question are known:

\begin{enumerate}
\item (Arvanitakis~\cite{Arvanitakis}) If $K$ is an almost totally
disconnected quasi Radon-Nikod\'ym compact, then $K$ is
Radon-Nikod\'ym compact.\\

\item (Avil\'es~\cite{cardinalb}) If $K$ is a quasi
Radon-Nikod\'ym compact of weight less than $\mathbf{b}$, then $K$
is Radon-Nikod\'ym compact.\\

\end{enumerate}

In addition, we also mention that some conditions are given by
Matou\v{s}kov\'a and Stegall \cite{MatSte} for a union of two
Radon-Nikod\'ym compacta to be Radon-Nikod\'ym compact. In this
paper we are mainly concerned with Arvanitakis' result, which
generalizes previous work of \cite{StegallRN} and
\cite{FabHeiMat}. We recall the concept of almost totally
disconnected compact. We denote by $\Sigma[0,1]^{\Gamma}$ the set
of all elements of $[0,1]^\Gamma$ with countable support and by
$\Sigma_0^1 [0,1]^{\Gamma}$ the subspace of $[0,1]^{\Gamma}$
formed by the elements such that all but
countably many coordinates belong to $\{0,1\}$.\\

\begin{defn}\label{definition of qtd}
A compact space is said to be almost totally disconnected if it
homeomorphic to a subset of
$\Sigma_0^1 [0,1]^{\Gamma}$ for some set $\Gamma$.\\
\end{defn}

This class contains both the classes of Corson compacta (compact
subsets of $\Sigma[0,1]^\Gamma$) and of totally disconnected
compacta. On the other hand, an idea of the limitations of
Arvanitakis' theorem may be suggested by the following remark,
which shows that the class of almost totally disconnected spaces
does not provide anything new with respect to Corson compacta,
when we restrict our attention to path-connected compacta, and
this includes the important case of probability measure spaces. We noticed
this fact in conversation with Ond\v{r}ej Kalenda during his visit to Murcia in March 2006.\\

\begin{prop}\label{pathconnected}
Every path-connected almost totally disconnected compact space is
a Corson compact.\\
\end{prop}

PROOF: Take $K\subset\Sigma_0^1 [0,1]^{\Gamma}$. We fix a point
$x\in K$ and without loss of generality we shall suppose that
$x_\gamma=0$ for all but countably many $\gamma$'s. In other
words, there exists a countable set $\Gamma_x\subset\Gamma$ such
that $x_\gamma = 0$ for every $\gamma\in\Gamma\setminus\Gamma_x$.
Now we take any other point $y\in K$ and we shall check that $y$
has also countable support. Since $K$ is path connected, there is
a separable and connected compact $L\subset K$ with $x,y\in
L\subset K$. Let $Q$ be a countable dense subset of $L$. For every
$q\in Q$ there is a countable set $\Gamma_q\subset\Gamma$ such
that $q_\gamma\in\{0,1\}$ for every $\gamma\in
\Gamma\setminus\Gamma_q$. The set $\Gamma_L = \bigcup_{q\in
Q}\Gamma_q$ is a countable subset of $\Gamma$ such that
$p_\gamma\in\{0,1\}$ for every $p\in L$ and every
$\gamma\in\Gamma\setminus\Gamma_L$. Since $L$ is connected, the
set $\{p_\gamma : p\in L\}$ must be connected for every
$\gamma\in\Gamma$. If we take
$\gamma\in\Gamma\setminus(\Gamma_L\cup\Gamma_x)$ then
$\{0\}\subset \{p_\gamma : p\in L\}\subset \{0,1\}$, so
connectedness implies $\{p_\gamma : p\in L\}=\{0\}$. Applying this
in particular to $p=y$, we found that $y_\gamma=0$ whenever
$\gamma\in\Gamma\setminus(\Gamma_L\cup\Gamma_x)$, so $y$ has
countable support.$\qed$\\

Apparently, we used a weaker hypothesis than path-connected in
this result, namely that every two points are contained in a
separable connected compact. However this is equivalent in this
context, because a separable connected compact which is almost
totally disconnected must be metrizable: Take
$K\subset\Sigma_0^1[0,1]$ and suppose that $D\subset K$ is a
countable dense subset of $K$. Then, we can find a countable
subset $\Gamma'\subset\Gamma$ such that $D\subset
[0,1]^{\Gamma'}\times\{0,1\}^{\Gamma\setminus\Gamma'}$. Hence
$K\subset [0,1]^{\Gamma'}\times\{0,1\}^{\Gamma\setminus\Gamma'}$
and the connectedness of $K$ implies that we have an embedding
$K\hookrightarrow [0,1]^{\Gamma'}$, so $K$ is metrizable. On the
other hand, the assumption of being path-connected cannot be
weakened to just being connected: several examples of connected
almost totally disconnected compacta which are not Corson will be described below.\\

The following Theorem~\ref{order fragmented are qtd} is the main
result of this note. Its proof is presented in Section~\ref{mainsection}.\\

\begin{thm}\label{order fragmented are qtd}
Let $K$ be a linearly ordered fragmentable compact. Then $K$ is almost totally disconnected.\\
\end{thm}

\begin{cor}\label{linearly ordered}
Let $K$ be a linearly ordered quasi Radon-Nikod\'ym compact. Then $K$ is a Radon-Nikod\'ym compact.\\
\end{cor}

A typical example of a linearly ordered compact which is not
fragmentable is the split interval (also known as double-arrow
space), that is, the set $K=[0,1]\times\{0,1\}$ ordered
lexicographically. Indeed any variant of the split interval on
which uncountably many points are splitted fails to be
fragmentable. The reason is that if $d$ is any quasi metric on
$K$, then there is an uncountable set $A\subset(0,1)$ and
$\varepsilon>0$ such that $d((x,0),(x,1))>\varepsilon$ for every
$x\in A$. If $B$ is a subset of $A$ in which every point of $B$ is
the limit of elements of $B$ both from the right and from the
left, then the set $B\times \{0,1\}$ fails to contain any relative
open subset of diameter less than
$\varepsilon$.\\

The class of linearly ordered compacta is a rather restrictive
class of compact spaces. For example, it is a result of Efimov and
\v{C}ertanov~\cite{EfiCerCorson}, with an alternative proof due to
Gruenhage \cite{GruenhageW}, that every linearly ordered Corson
compact space is metrizable. In the view of this result and also
of Proposition~\ref{pathconnected}, one may be suspicious about
real application of Theorem~\ref{order fragmented are qtd}. This
is not the case and indeed one of the examples of Radon-Nikod\'ym
compact proposed by Namioka~\cite{Namioka} is the so-called
extended long line. This is a linearly ordered compact obtained
from the ordinals less or equal to $\omega_1$ by inserting a copy
of the interval $(0,1)$ between every two consecutive countable
ordinals. More examples of linearly ordered Radon-Nikod\'ym
compacta are constructed in Section~\ref{secondsection}, where
Corollary~\ref{linearly
ordered} will find application.\\

It is an open question for us whether every fragmentable linearly
ordered compact must be a Radon-Nikod\'ym compact.\\

\section{Proof of the main theorem}\label{mainsection}

We begin with a couple of lemmas, stating reformulations of the
concept of almost totally disconnected compact, the second of
them in the framework of linearly ordered compacta.\\

\begin{lem}\label{characterization qtd}
For a compact space $K$ the following are equivalent:
\begin{enumerate}
\item $K$ is almost totally disconnected. \item There is a
collection $\{(F_{i},H_{i})\}_{i\in I}$ of pairs of closed subsets
of $K$ such that
\begin{enumerate}
\item $F_{i}\cap H_{i}=\emptyset$ for all $i\in I$. \item For any
$x$ in $K$, the set $\{i\in I : x\not\in F_{i}\cup H_{i}\}$ is
countable. \item For any two different points $x,y$ in $K$ there
is some $i\in I$ such
that $x\in F_{i}$ and $y\in H_{i}$ or viceversa.\\
\end{enumerate}
\end{enumerate}
\end{lem}

PROOF: For $(1)\Rightarrow (2)$, suppose $K\subset \Sigma_0^1
[0,1]^{\Gamma}$. For each $\gamma\in\Gamma$ and each pair $r,s$ of
rational numbers with $0\leq r<s\leq 1$, call $i=(\gamma,r,s)$,
$$F_{i} = \{ x\in K : x_{\gamma}\leq r\}$$
$$H_{i} = \{x\in K : x_{\gamma}\geq s\}$$
These $(F_{i},H_{i})$ satisfy all desired conditions.\\

Conversely, suppose we are given a family like in (2). For each
$i$, by Tietze's theorem, there is a continuous map $f_{i}:K\To
[0,1]$ such that $f_{i}(F_{i})=\{0\}$ and $f_{i}(H_{i}) = \{1\}$.
In this case we have an embedding $f:K\To
\Sigma_0^1[0,1]^{I}$ given by $f(x) = (f_{i}(x))_{i\in I}$.\hfill$\square$\\

\begin{lem}\label{QTD order}
Let $(K,\leq)$ be a linearly ordered compact space. The following
are equivalent:
\begin{enumerate}
\item K is almost totally disconnected. \item There is a
collection $\{(a_{i},b_{i})\}_{i\in I}\subset K\times K$ such that
\begin{enumerate}
\item $a_{i} < b_{i}$ for every $i\in I$. \item For all $x$ in
$K$, the set $\{i\in I : a_{i}<x<b_{i}\}$ is countable. \item For
all $x<y$ in $K$, there is some $i\in I$ such that $x\leq
a_{i}<b_{i}\leq
y$.\\
\end{enumerate}
\end{enumerate}
\end{lem}
PROOF: Clearly $(2)$ implies $(1)$ because $F_{i} =
]-\infty,a_{i}]$ and $H_{i}=[b_{i},+\infty[$ satisfy the
conditions of Lemma~\ref{characterization qtd}. Conversely,
suppose that we have a family $(F_{j},H_{j})_{j\in J}$ of couples
of closed subsets of $K$ satisfying the conditions of
Lemma~\ref{characterization qtd}. Take as $\{(a_{i},b_{i})\}_{i\in
I}$ the set of all pairs in $K\times K$ such that\begin{enumerate}
\item $a_{i}<b_{i}$ \item There is some $j(i)\in J$ such
that\begin{enumerate}
 \item Either $a_{i}\in F_{j(i)}$ and $b_{i}\in H_{j(i)}$, or viceversa, $b_{i}\in F_{j(i)}$ and $a_{i}\in H_{j(i)}$.
 \item There is no $x\in F_{j(i)}\cup H_{j(i)}$ such that
 $a_{i}<x<b_{i}$.\\
 \end{enumerate}
\end{enumerate}

First, we check that for every $x\in K$, the set $\{i\in I :
a_i<x<b_i\}$ is countable. Notice that whenever $a_i<x<b_i$, then
$x\not\in F_{j(i)}\cup H_{j(i)}$ and we know that, since the
family $\{(F_j,H_j)\}_{j\in J}$ satisfies condition (b) of
Lemma~\ref{characterization qtd}, the set $\{j\in J : x\not\in
F_j\cup H_j\}$ is countable. The fact that $\{i\in I :
a_i<x<b_i\}$ is countable follows now from the observation that
whenever $j(i) = j(i')$ and $i\neq i'$, the intervals $]a_i,b_i[$
and
$]a_{i'},b_{i'}[$ are disjoint.\\

Second, we check condition (c) of the lemma. Take $x<y$. Since
condition (c) of Lemma~\ref{characterization qtd} is satisfied, we
suppose that there is some $j\in J$ such that $x\in F_j$ and $y\in
H_j$. Let $z=\max\{t\in F_j : t\leq y\}$ and $z'=\min\{t\in H_j :
z\leq t\}$. Then,
$(z,z')$ equals some $(a_i,b_i)$ and $x\leq a_i<b_i\leq y$.\hfill$\square$\\

We pass now to the proof of Theorem~\ref{order fragmented are qtd}
itself. Let $K$ be a linearly ordered compact and let $d$ be a
quasi metric which fragments $K$. We construct our family
$\{(a_{i},b_{i})\}\subset K\times K$ as in Lemma~\ref{QTD order}
as follows. First, let $\{(a_{i},b_{i})\}_{i\in I_{0}}$ be the set
of all pairs of immediate successors (that is, all $a_i<b_i$ such
that the open interval $]a_{i},b_{i}[$ is empty). For $n\geq 1$,
by virtue of Zorn's Lemma, we can choose a family
$(a_{i},b_{i})_{i\in I_{n}}$
which is maximal for the following properties:\\
\begin{enumerate}
\item $a_i<b_i$ for every $i\in I_n$. \item The $d$-diameter of
the open interval $]a_{i},b_{i}[$ is less than $1/n$. \item
$[a_{i},b_{i}]\cap [a_{j},b_{j}] = \emptyset$ whenever $i,j$ are
different indices in $I_{n}$.\\
\end{enumerate}

We take $I=\bigcup_{n=0}^{\infty}I_{n}$ and $(a_{i},b_{i})_{i\in
I}$ as the family required in Lemma~\ref{QTD order}. Condition (a)
of Lemma~\ref{QTD order} is clearly satisfied and condition (b)
follows from property (3) in the definition of $I_{n}$. Only
condition (c) needs
to be checked. Take $x<y$, and we  suppose that\\

(A)  there is no index $i\in I$ such that $x\leq a_{i}<b_{i}\leq y$.\\

This implies, by the definition of $I_0$, that no immediate
successors can occur between $x$ and
$y$. This means that the interval $[x,y]$ is connected (and so all its subintervals).\\

It is not possible that for all $n$, there is $j\in I_{n}$ such
that $]x,y[\subset ]a_{j},b_{j}[$. This is because, by property
(2) of $I_{n}$, the $d$-diameter of
$]x,y[$ would be 0, which is a contradiction. Hence,\\

(B)  for some fixed $n_0\in\omega$, there is no $j\in I_{n_0}$ such that $]x,y[\subset ]a_{j},b_{j}[$.\\

Claim 1: There exists $i_0\in I_{n_0}$ such that either $x<a_{i_0}<y$ or $x<b_{i_0}<y$.\\

Proof of the claim: If such an $i_0$ does not exist, assertion (B)
implies that $]x,y[\cap]a_i,b_i[=\emptyset$ for all $i\in
I_{n_0}$. Since $d$ fragments $K$ there is a nonempty open
interval $]u,v[\subset]x,y[$ of $d$-diameter less than
$\frac{1}{n_0}$. By passing to a subinterval (recall that all
intervals are connected now) it can be supposed that even
$[u,v]\subset]x,y[$ and then, the pair $(u,v)$ could be added to
the family $\{(a_i,b_i)\}_{i\in I_{n_0}}$ in contradiction with
its maximality.\\

Without loss of generality, we assume that $x<a_{i_0}<y$ in Claim~1.\\

Claim 2: There exists $j_0\in I_{n_0}$ such that $x<b_{j_0}<a_{i_0}<y$.\\

Proof of the claim: Again, if such a $j_0$ does not exist, then
$]x,a_{i_0}[\cap ]a_j,b_j[$ is empty for all $j\in I_{n_0}$ and,
because of the fragmentability condition, we can find an interval
$[u,v]\subset ]x,a_{i_0}[$ of $d$-diameter less than
$\frac{1}{n_0}$. In this case, the couple $(u,v)$ could be added
to the family
$\{(a_j,b_j)\}_{j\in I_{n_0}}$ in contradiction with its maximality.\\

Claim 3: There exists $k_0\in I_{n_0}$ such that
$x<b_{j_0}<a_{k_0}<b_{k_0}<a_{i_0}<y$.\\

Proof of the claim: Similarly, if such a $k_0$ did not exist, then
$]b_{j_0},a_{i_0}[\cap ]a_k,b_k[=\emptyset$ for all $k\in
I_{n_0}$. Since $d$ fragments $K$, there should be an interval
$[u,v]\subset ]b_{j_0},a_{i_0}[$ of $d$-diameter less than
$\frac{1}{n_0}$ which leads once again to a contradiction.\\

Claim 3 is incompatible with assumption (A) and this finishes the
proof of Theorem~\ref{order fragmented are qtd}.\hfill$\square$\\

Notice that we did not use the full strength of the definition of
fragmentability. We just needed that $d(x,y)>0$ if $x\neq y$ and
that every \emph{interval} contains an open subinterval of
$d$-diameter less than $\varepsilon$ for every $\varepsilon>0$.

\section{Examples of linearly ordered Radon-Nikod\'ym
compacta}\label{secondsection}

We recall a different characterization of quasi Radon-Nikod\'ym
compacta. A metric $d:K\times K\To [0,+\infty)$ on the compact
space $K$ is called a Reznichenko metric \cite[p.
104]{ArkhangelskiiGT2} if for every two different points $x,y\in
K$ there exist neighborhoods $U$ and $V$ of $x$ and $y$
respectively such that $\inf\{d(u,v) : u\in U, v\in V\}>0$. The
following theorem is due to Namioka~\cite{Namiokanote}:\\

\begin{thm}\label{Namiokanoteresult}
A compact space $K$ is quasi Radon-Nikod\'ym if and only if there
exists a Reznichenko metric on $K$ which fragments $K$.\\
\end{thm}

We present now a method for constructing linearly ordered
Radon-Nikod\'ym compact spaces inspired on Ribarska's
characterization of fragmentability~\cite{Ribarska}. We consider
$\{T_n : n=1,2,\ldots\}$ to be a sequence of well ordered sets
such that $T_n\subset T_{n+1}$. Without loss of generality, we
shall assume that all $T_n$'s have the same minimum and the same
maximum. Let $T$ be the linearly ordered set
$T=\bigcup_{n=1}^\infty T_n$ and let $\bar{T}$ be the completion
of $T$ (by the completion of $T$ we mean the only linearly ordered
set $\bar{T}$ such that $T\subset \bar{T}$, $\bar{T}$ is compact
in the order topology and $]x,y]\cap T\neq\emptyset$ for every
$x,y\in \bar{T}$, $x<y$). Then, $\bar{T}$ is a linearly ordered
compact
space and moreover:\\

\begin{thm}\label{espacio T}
The space $\bar{T}$ is Radon-Nikodým compact.\\
\end{thm}

This theorem produces different examples of linearly ordered
Radon-Nikod\'ym compacta depending on the growing sequence of well
ordered sets $T_1\subset T_2\subset\cdots$ that we may take as a
basis. Before passing to the proof, we shall have a look at how
these different constructions may look like. Notice that $\bar{T}$
is connected whenever for all $x<y$ in $T_n$ there exists $z\in T_{n+1}$ such that $x<z<y$.\\

\begin{itemize}

\item If $T_0 = \{0,1\}$ and $T_{n+1}$ is constructed by adding a
single new point between every two consecutive elements of $T_n$,
then $\bar{T} = [0,1]$.\\

\item If $T_0$ is the set of all ordinals which are less than or
equal to $\omega_1$ and again $T_{n+1}$ is constructed by adding a
single new point between every two consecutive elements of $T_n$,
then $\bar{T}$ is the extended long line.\\

\item If $T_0$ is the set of all ordinals which are less than or
equal to $\omega_1$ and $T_{n+1}$ is constructed by adding a copy
of the set of all countable ordinals between every two consecutive
elements of $T_n$,
then $\bar{T}$ has no metrizable open subsets, since every open interval contains a copy of $\omega_1$.\\

\end{itemize}

Proof of Theorem~\ref{espacio T}: For $x,y\in \bar{T}$, $x<y$ we define:\\
$$d(x,y)=\frac{1}{\min\{n: ]x,y]\cap T_n\neq\emptyset\}}$$
and also $d(x,y)=0$ if $x=y$, and $d(x,y)=d(y,x)$ if $x>y$. Observe that: \\

\begin{enumerate}
\item Since $\bar{T}$ is the completion of $T$, if $x<y$ then
$T\cap ]x,y]$ is nonempty. This implies that $d(x,y)$ exists and
is a positive real whenever $x\neq
y$.\\

\item An easy case-by-case consideration proves that $d$ verifies
triangle inequality: $d(x,z)\leq d(x,y)+d(y,z)$.\\

\item The metric $d$ is a Reznichenko metric, that is, every two
different points have neighborhoods at a positive $d$-distance.
Namely, if $x<y$ and there is some $z\in]x,y[$, then $]x,z]\cap T$
is nonempty and there is $u\in T_n\cap]x,z]$ for some $n$. In this
case $]-\infty,u[$ and $]u,+\infty[$ are neighborhoods at at least
$\frac{1}{n}$-$d$-distance. The other possibility is that $]x,y[$
is empty. Then $y\in T$ and hence, $y\in T_n$ for some $n$. In
this case $(-\infty,x]$ and $[y,+\infty[$ are neighborhoods of $x$
and $y$ at at
least $\frac{1}{n}$-$d$-distance.\\

\item The metric $d$ fragments $\bar{T}$. Given $L$ a closed
subset of $\bar{T}$ of more than one point, and $n\in\omega$, let
$x=\min(L)$ and $y=\min\{ z\in T_n : z>x\}$. Then, $L\cap [x,y[$
is a nonempty relative open subset of
$L$ of $d$-diameter less than $\frac{1}{n}$.\\

\end{enumerate}

It follows from Theorem~\ref{Namiokanoteresult} that $K$ is quasi
Radon-Nikod\'ym compact and hence, by Corollary~\ref{linearly
ordered}, it is Radon-Nikodým compact.\hfill$\square$\\

We point out that the metric $d$ that we have defined is not lower
semicontinuous except in trivial cases and hence, the use of
Corollary~\ref{linearly ordered} cannot be avoided unless we want
to embark in the construction of a more complicated metric. In
some cases, a nice lower semicontinuous fragmenting metric may be
available, as for the extended long line~\cite{Namioka}.

\bibliographystyle{amsplain}

\end{document}